 \documentclass[final,5p,times,twocolumn,fleqn]{elsarticle}

\usepackage{amssymb}
\usepackage[cmintegrals]{newtxmath}
\usepackage{bm}
\usepackage{amsmath}
\usepackage{graphicx}
\usepackage{graphics}
\usepackage{color}
\usepackage[hyperindex]{hyperref}

\newtheorem{theorem}{Theorem}[section]

\newtheorem{definition}[theorem]{Definition}
\newtheorem{problem}[theorem]{Problem}

\newtheorem{remark}[theorem]{Remark}

\newtheorem{assumption}[theorem]{Assumption}

\newtheorem{lemma}[theorem]{Lemma}

\journal{}
\begin{document}
\newcommand{\phase}{\mbox{$\mathbb{\delta}$}}
\newcommand{\realpart}{\mbox{\rm Re}}
\newcommand{\da}{\mbox{$\rm DA$}}
\newcommand{\domainofattraction}{\mbox{$DA$}}
\begin{frontmatter}
%
%

\title{Transient Performance of Power Systems\\ with Distributed Power-Imbalance 
Allocation Control}
\tnotetext[titlenote]{This work is supported in part by the research project of  the Fundamental Research Funds of Shandong University under Grant 2018HW028 and in part by the 
Foundation for Innovative Research Groups of National Natural Science Foundation of China under Grant 61821004.}

%
\author[1]{Kaihua Xi\corref{cor1}}
\cortext[cor1]{Corresponding author}
\ead{kxi@sdu.edu.cn}
\address[1]{School of Mathematics, Shandong University, Jinan, 250100, Shandong, China}
\author[2]{Hai Xiang Lin}
\ead{H.X.Lin@tudelft.nl}
\author[2]{Jan H. van Schuppen}
\ead{jan.h.van.schuppen@tudelft.nl}
\address[2]{Delft Institute of Applied 
Mathematics, Delft University of Technology, Van Mourik Broekmanweg 6,
2628 XE Delft, The Netherlands.}

\begin{abstract}
We investigate the sensitivity of the transient performance of power systems controlled by \emph{Distributed Power-Imbalance Allocation Control} (DPIAC) on the coefficients of the 
control law. We measure the transient performance of the frequency deviation and the control cost by the $\mathcal{H}_2$ norm. 
Analytic formulas are derived for the $\mathcal{H}_2$ norm of the transient performance of a power system with homogeneous parameters and with
a communication network of the same topology as the power network. It is shown that the transient performance of the frequency can be greatly improved by accelerating the convergence to the optimal steady state through the control gain coefficients, which however requires a higher control 
cost. Hence, in DPIAC, there is a trade-off between the frequency deviation and the control cost which is determined by the control gain coefficients. In addition, by increasing one of the control gain coefficients, the behavior of the state approaches that of a centralized control law. These analytical results are validated through a numerical simulation of the IEEE 39-bus system in which the system parameters are heterogeneous.
\end{abstract}
%
%
%
\end{frontmatter}
%
%
%
\section{Introduction}
The power system is expected to keep the 
frequency within a small range around the nominal value so as to avoid damages to electrical devices. This is accomplished by regulating the active power injection of sources.  
Three forms of frequency control can be distinguished from fast to slow timescales, i.e., primary control, secondary control, 
and tertiary control
\cite{Kundur1994,Wood}. 
Primary frequency control has a control objective 
to maintain the synchronization of the frequency based on local feedback at each power generator. However, the synchronized frequency of the entire power system with primary controllers may still deviate from its nominal value. Secondary frequency control restores the synchronized frequency to its nominal value and is operated on a slower time scale than primary control. Based on the predicted power demand, tertiary control determines the set points for both
primary and secondary control over a longer period than used in secondary control. The operating point is usually the solution of an optimal power flow problem.


The focus of this paper is on the secondary frequency which has been traditionally 
actuated by the passivity-based method \emph{Automatic Genertaion Control}~(AGC) for half a century. Recently, considering the on-line economic power dispatch in the secondary frequency control between all the controllers in the power system \cite{Wood}, 
various control methods are proposed for the secondary frequency control. 
These include passivity-based centralized control methods such as \emph{Gather-Broadcast} Control \cite{Dorfler2016}, and distributed control methods 
such as the \emph{Distributed Average Integral} (DAI)\cite{Zhao2015}, primal-dual algorithm based distributed method
such as the \emph{Economic Automatic Generation Control} (EAGC) \cite{EAGC}, \emph{Unified Control} \cite{UC} \emph{etc..} 
However, the primary design objectives of these methods focus on the steady state only. As investigated in our previous study in \cite{PIAC3}, the 
corresponding closed-loop system may have a poor transient performance even though the control objective of reaching the steady state is achieved, e.g. \cite{quadratic_performance,PIAC3}. For example, from the global perspective of the entire power system, the passivity based methods, e.g., AGC, GB and DAI, are actually a form of integral control. A drawback of integral control 
is that large integral-gain coefficients may result in extra oscillations due to the \emph{overshoot of the control input} while small gain coefficients result in a slow convergence speed towards a steady state. For 
instance, an overshoot problem occurred after the blackout which happened in UK in August 2019~\cite{England2019}.
\par
The continuous increase of integration of renewable energy into the power systems, which may bring serious fluctuations, asks for more attentions on the transient performance of the power systems.
For the secondary frequency control, the way to 
improve the transient performance of the traditional methods is
to tune the control gain coefficients either by obtaining satisfactory 
eigenvalues of the linearized closed-loop system or by using a control law 
based on $\mathcal{H}_2$ or $\mathcal{H}_{\infty}$ control synthesis      \cite{Hassan_Bevrani,Hinfinty_control_overshoot}. However, besides the complicated computations, the improvement of the 
transient performance is still limited because it also depends on the structure
of the control laws. In order to improve the transient performance, sliding-mode-based control laws, e.g.,\cite{DeSlide,VRDOLJAK2010514} and fuzzy control-based control laws,e.g., \cite{CAM2005233} are proposed, which are all able to shorten the transient phase without the overshoot. However,
those control laws use either centralized or decentralized control structure without considering economic power dispatch. 
Concerning the transient performance and the balance of the advantages of the centralized and distributed control structure, the authors have proposed a multilevel control method, \emph{Multi-Level Power-Imbalance Allocation Control} (MLPIAC) \cite{MLPIAC} for the secondary 
frequency control, which is suitable for large scale power systems.  There are two special cases of MLPIAC, a centralized control called \emph{Gather-Broadcast Power-Imbalance Allocation Control} (GBPIAC), 
and a distributed control called \emph{Distributed Power-Imbalance Allocation Control}~(DPIAC). Numerical study with comparison to the existing methods show that the overshoot problem can be avoided by MLPIAC, thus the transient performance can be improved by accelerating the 
convergence of the state~\cite{PIAC3,MLPIAC}. However, the analysis is incomplete, which lacks of quantifying the 
impact of these control coefficients on the transient performance. 
\par 
The $\mathcal{H}_2$ norm of a time-invariant linear input-output system reflects the response of the output to the input, which has been
widely used to study the response of power systems to disturbances, e.g.,
the performance analysis of secondary frequency control methods in \cite{DAPI_performance,quadratic_performance,Wu2016,Dorfler2019}, the optimal virtual inertia placement in Micro-Grids \cite{optimal_inertia_placement}, and the cyber network design for secondary frequency control \cite{LGUO}.
 
In this paper, we focus on the distributed control law DPIAC, and analyze the impact of the control coefficients on the transient performance of the frequency deviation and the control cost after a disturbance. For comparison with DPIAC, we also investigate the transient performance of the centralized control method GBPIAC. We measure the transient performance by the $\mathcal{H}_2$ norm. We will show analytically and numerically that (1) the transient performance
can be improved by tuning the control gain coefficients monotonically; (2) there is a trade-off between the transient performance of frequency deviation 
and the control cost, which is determined by the control coefficients; and (3)
the performance of the distributed control approaches
to that of the centralized control as a gain coefficient is increased. 
The main contributions of this paper are,
\begin{enumerate}[(i)]
\item analytic formulas for how the transient performance of the 
frequency deviation and the control cost depends on the control gain coefficients; 
\item a numerical study of  the transient performance and its dependence
on the control gain coefficients. 
\end{enumerate}
\par 
The paper is organized as follows. We first introduce the model
model of the power system, GBPIAC and DPIAC in section \ref{chapter5:Sec:model},
then formulate the problem of this paper with introduction 
of the $\mathcal{H}_2$ norm in section \ref{chapter5:Sec:problem}. 
We calculate the corresponding $\mathcal{H}_2$ norms and analyze the impact of the control coefficients on the transient performance of the frequency deviation and the control cost in section \ref{chapter5:Sec:coherence} and verify 
the analysis by simulations in section \ref{section:simulationstudy}. Finally,
we conclude with remarks in section \ref{chapter5:Sec:conclusion}. 

\section{The secondary frequency control laws}\label{chapter5:Sec:model}

The transmission network of a power system can be described by a graph $\mathcal{G}=(\mathcal{V},\mathcal{E})$ with nodes $\mathcal{V}$ and edges
$\mathcal{E}\subseteq \mathcal{V}\times \mathcal{V}$,  where a node represents a bus and edge $(i,j)$ represents the direct transmission line between node $i$ and $j$. The buses can connect to synchronous machines, frequency dependent power sources (or loads), or passive loads. We focus on
a power system with loss-less transmission lines and denote the susceptance of the transmission line by 
 $\hat{B}_{ij}$ for $(i,j)\in\mathcal{E}$. The set of the buses of the synchronous machines, of frequency dependent power sources, of passive 
 loads are denoted by $\mathcal{V}_M,\mathcal{V}_F,\mathcal{V}_P$ respectively, thus $\mathcal{V}=\mathcal{V}_M\cup\mathcal{V}_F\cup\mathcal{V}_P$.
 The dynamics of the system can be modelled by the following \emph{Differential Algebraic Equations} (DAEs), e.g., \cite{Dorfler2016,PIAC3}, 
\begin{subequations}\label{chapterModel:eq:system1}
\begin{align}
  \hspace{-10pt}
  \dot{\theta}_i&=\omega_i,~i\in\mathcal{V}_M\cup\mathcal{V}_F,\\
    \hspace{-10pt}
M_{i}\dot{\omega}_{i}&=P_{i}-D_{i}{\omega}_{i}-\sum_{j\in \mathcal V}{K_{ij}\sin{(\theta_{i}-\theta_{j})}}+u_{i}, 
i\in \mathcal{V}_M,\label{eq:syn}\\
  \hspace{-10pt}
0&=P_{i}-D_{i}\omega_{i}-\sum_{j\in \mathcal V}{K_{ij}\sin{(\theta_{i}-\theta_{j})}}+u_{i}, i\in \mathcal{V}_F,\label{eq:frq}\\
  \hspace{-10pt}
0&=P_{i}-\sum_{j\in \mathcal V}{K_{ij}\sin{(\theta_{i}-\theta_{j})}}, i\in \mathcal{V}_P, \label{eq:pass} 
 \end{align}
\end{subequations}
where $\theta_i$ is the phase angle at node $i$, $\omega_i$ is the 
frequency deviation from the nominal value (e.g., 50 or 60 Hz), $M_i>0$ is the moment of inertia of the synchronous machine, 
 $D_i>0$ is the droop control coefficient, $P_i$ is the 
power supply (or demand), $K_{ij}=\hat{B}_{ij}V_iV_j$ is the effective 
susceptance of the transmission line between node $i$ and $j$, $V_i$ is the voltage, $u_i$ is the input for the secondary control. 
The nodes in $\mathcal{V}_M$ and $\mathcal{V}_F$ are assumed to be equipped with secondary frequency controllers, denoted by $\mathcal{V}_K=\mathcal{V}_M\cup\mathcal{V}_F$. 
Since the control of the voltage and the frequency can be decoupled when the transmission lines are lossless \cite{Simpson-Porco2016}, we do not model the dynamics
of the voltages and assume 
the voltages are constant which can be derived from a power flow calculation \cite{Kundur1994}.

The synchronized frequency deviation can be expressed as 
\begin{eqnarray}
 \omega_{syn}=\frac{\sum_{i\in \mathcal{V}}{P_i}+\sum_{i\in\mathcal{V}_K}{u_i}}{\sum_{i\in \mathcal{V}_M \cup\mathcal{V}_F }{D_i}}. 
 \label{chapterModel:eq:sync}
\end{eqnarray}
The condition $\omega_{syn}=0$ at a steady state can be satisfied by solving the economic power dispatch problem in the secondary frequency control \cite{Wood},   
 \begin{eqnarray}
 &&\min_{u_i\in R} \sum_{i\in\mathcal{V}_K }J_i(u_i)
 \label{eq:optimal3}\\
 \nonumber
 &&s.t. ~~\sum_{i\in\mathcal{V}}P_i+\sum_{i\in\mathcal{V}_K }{u_i}=0,
 \end{eqnarray}
where
$J_i(u_i)=\frac{1}{2}\alpha_i u_i^2$ represents the control cost 
at node $i$. The constraints on the control input at each node is not 
included into (\ref{eq:optimal3}) as in \cite{Wood,MLPIAC} based on the assumption that the constraint of the capacity is not triggered by the disturbances 
in the time-scale of the secondary frequency control. This assumption allows the use of the $\mathcal{H}_2$ norm in the transient performance analysis.

A necessary condition for the solution of the optimization problem (\ref{eq:optimal3}) is that 
the \emph{marginal costs} $dJ_i(u_i)/du_i$ of the nodes are all identical, i.e.,
\begin{eqnarray*}
 \alpha_i u_i=\alpha_j u_j,~ \forall~i,~j\in\mathcal{V}_K.
\end{eqnarray*}

For a secondary control law with the objective of (\ref{eq:optimal3}), 
it is required that the total control input $u_s(t)=\sum_{i\in\mathcal{V}_K }{u_i}$ converges to the unknown $-P_s=-\sum_{i\in\mathcal{V}}P_i$ 
and the marginal costs achieve a consensus at the steady state. 
For the overview of the secondary frequency control laws, see \cite{overviewSecondaryFreq}. To obtain a good transient performance, a fast convergence of the control inputs to the optimal solution of (\ref{eq:optimal3}) is critical, which may introduce the overshoot problem. To avoid the overshoot problem, 
MLPIAC has been proposed in \cite{MLPIAC} with two special cases, GBPIAC and DPIAC.  In this paper, we focus on the impact of the control gain coefficients of DPIAC on the transient performance and compare with that of GBPIAC. The following assumption on the connectivity of the communication network is required to realize the coordination control. 
\begin{assumption}
 For the power system (\ref{chapterModel:eq:system1}), there exists a undirected  communication network such that all the nodes in $\mathcal{V}_K$ are connected. 
\end{assumption}

%
%

The definition of the distributed
method DPIAC and the centralized method GBPIAC follow. 

\begin{definition}[GBPIAC]
Consider the power system (\ref{chapterModel:eq:system1}), the 
\emph{Gather-Broadcast Power-Imbalance Allocation Control}
(GBPIAC) law is defined as \cite{MLPIAC}
\begin{subequations}\label{eq:GBPIAC}
 \begin{align}
  \dot{\eta}_s&=\sum_{i\in\mathcal{V}_M\cup\mathcal{V}_F}D_i\omega_i,\label{eq:controlPIAC3a}\\
  \dot{\xi}_s&=-k_1(\sum_{i\in\mathcal{V}_M}M_i\omega_i+\eta_s)-k_2\xi_s,\label{eq:controlPIAC3b}\\
   u_i&=\frac{\alpha_s}{\alpha_i}k_2\xi_s,~i\in\mathcal{V}_K,\label{eq:uxi}
 \end{align}
\end{subequations}
where $\eta_s\in\mathbb{R},~\xi_s\in\mathbb{R}$ are state variables of the central controller,  $k_1, k_2$ are positive control gain coefficients, $\alpha_i$ is the 
control price at node $i$ as defined in the optimization problem (\ref{eq:optimal3}), $\alpha_s=(\sum_{i\in\mathcal{V}_K}1/\alpha_i)^{-1}$
is a constant.
\end{definition}



\begin{definition}[DPIAC]\label{def:DPIAC}
Consider the power system (\ref{chapterModel:eq:system1}), define the 
\emph{Distributed Power-Imbalance Allocation Control} (DPIAC)  law as,
\begin{subequations}\label{eq:DPIAC}
 \begin{align}
\dot{\eta}_i&=D_i\omega_i+k_3\sum_{j\in\mathcal{V}_K}{l_{ij}{(k_2\alpha_i\xi_i-k_2\alpha_j\xi_j)}},\label{eq:DPIAC1}\\
\dot{\xi}_i&=-k_1(M_i\omega_i+\eta_i)-k_2\xi_i,\label{eq:DPIAC2}\\
u_i&=k_2\xi_i, \label{eq:DPIAC2u}
 \end{align}
\end{subequations}
for node $i\in\mathcal{V}_K$, where
$\eta_i\in\mathbb{R}$ and $\xi_i\in\mathbb{R}$ are state variables of the local controller at node $i$, 
$k_1, k_2$ and $k_3$ are positive gain coefficients, $(l_{ij})$ defines a weighted undirected communication network with 
Laplacian matrix $(L_{ij})$
\begin{equation*}
 L_{ij}=\begin{cases}
         -l_{ij}, &i\neq j,\\
         \sum_{k\neq i}l_{ik}, &i=j,
        \end{cases}
\end{equation*}
and $l_{ij}\in[0,\infty)$ is the weight of the communication line connecting node $i$ and $j$. The marginal cost at 
node $i$ is $\alpha_i u_i=k_2\alpha_i\xi_i$. 
\end{definition}



Without the coordination on the marginal costs, DPIAC reduces to a decentralized control method as follows.
\begin{definition}[DecPIAC]
Consider the power system (\ref{chapterModel:eq:system1}), the \emph{Decentralized 
Power-Imbalance Allocation Control} (DecPIAC) 
is defined as,
 \begin{subequations}\label{decentralized}
 \begin{align}
\dot{\eta}_i&=D_i\omega_i,\\
\dot{\xi}_i&=-k_1(M_i\omega_i+\eta_i)-k_2\xi_i,\\
u_i&=k_2\xi_i,
 \end{align}
\end{subequations}
for node $i\in\mathcal{V}_K$, where $k_1$ and $k_2$ are positive gain coefficients. 
\end{definition}

For a fast recovery from an imbalance while avoiding the overshoot problem, the control gain coefficient $k_2$ should satisfy $k_2\geq 4k_1$. For details of the configuration of $k_1$ and $k_2$, we refer to \cite{MLPIAC}. In this paper, we set $k_2=4k_1$  in the following analysis so as to simplify the deduction of the explicit formula which shows 
the impact of the control coefficients on the transient performance.
For the control procedure and the asymptotic stability  of GBPIAC, DPIAC, see \cite{MLPIAC}. Fo the control law MLPIAC, see \cite{MLPIAC}. 
\section{Problem formulation}\label{chapter5:Sec:problem}

With $k_2=4k_1$, there are two control gain coefficients $k_1$ and $k_3$ in DPIAC.  
We focus on the following problem. 
\begin{problem}\label{chapter5:question}
How do the coefficients $k_1$  and $k_3$ influence the transient performance 
of the frequency deviation and control cost in the system (\ref{chapterModel:eq:system1}) controlled by DPIAC? 
\end{problem}
\par 
To address Problem \ref{chapter5:question}, we introduce the $\mathcal{H}_2$ norm
to measure the transient performance, which is defined as follows. 
\begin{definition}\label{definition_H2norm}
Consider a linear time-invariant system, 
\begin{subequations}\label{Appendix:generalform}
 \begin{align}
  \dot{\bm x}&=\bm A\bm x+\bm B\bm w,\\
  \bm y&=\bm C\bm x, 
 \end{align}
\end{subequations}
where $\bm x\in\mathbb{R}^n$, $\bm A\in\mathbb{R}^{n\times n}$ is Hurwitz, $\bm B\in\mathbb{R}^{n\times m}$,  $\bm C\in\mathbb{R}^{z\times n}$, 
the input is denoted by $\bm w\in\mathbb{R}^m$ and the output of the system is denoted by $\bm y\in\mathbb{R}^z$. 
The squared $\mathcal{H}_2$ norm of the transfer matrix $\bm G$ of the mapping $(\bm A,\bm B,\bm C)$ from the input $\bm w$ to the output $\bm y$ is defined as 
\begin{subequations}\label{Appendix:H2norm}
 \begin{align}
  &||\bm G||^2_2=\text{tr}(\bm B^T\bm Q_o\bm B)=\text{tr}(\bm C\bm Q_c\bm C^T),\\
  &\bm Q_o\bm A+\bm A^T\bm Q_o+\bm C^T\bm C=\bm 0,\label{Lyapunov:equation}\\
  &\bm A\bm Q_c+\bm Q_c\bm A^T+\bm B\bm B^T=\bm 0,
 \end{align}
\end{subequations}
where $\text{tr}(\bm \cdot)$ denotes the trace of a matrix, $\bm Q_o,\bm Q_c\in\mathbb{R}^{n\times n}$ are the \emph{observability Grammian} of $(\bm C,\bm A)$ and \emph{controllability 
Grammian} of $(\bm A,\bm B)$ respectively \cite{H2norm_another_form},\cite[chapter 2]{H2norm_book_toscano}. 
\end{definition}
\par 
The $\mathcal{H}_2$ norm can be interpreted as follows.
When the input $\bm w$ is modeled as the Gaussian white noise such that $w_i\sim N(0,1)$ for all $i=1,\cdots,m$ and for all $i\neq j$, the scalar Brownian motions $w_i$ and $w_j$ are independent, the matrix $\bm Q_{v}=\bm C\bm Q_c\bm C^T$ is the variance matrix of the output at the steady state \cite[Theorem 1.53]{karatzas}, i.e., 
\begin{eqnarray*}
 \bm Q_v=\lim_{t\rightarrow \infty}E[\bm y(t)\bm y^T(t)]
\end{eqnarray*}
where $E[\cdot]$ denotes the expectation. Thus 
\begin{eqnarray}
 \|\bm G\|^2_2=\text{tr}(\bm Q_v)=\lim_{t\rightarrow \infty}E[\bm y(t)^T\bm y(t)]. \label{chapter5:expected_value}
\end{eqnarray}
For other interpretations, see \cite{H2norm}.


\par 
There are so many parameters
which influence the transient performance of the system that it is hard to 
deduce an explicit formula of the $\mathcal{H}_2$ norm for the closed-loop system 
when the parameters are heterogeneous. To simplify the analysis and focus on the 
impact of the control gain coefficients, we make the following assumption. 
\begin{assumption}\label{assumption_a1}
For GBPIAC and DPIAC, assume that $\mathcal{V}_F=\emptyset$, $\mathcal{V}_P=\emptyset$ and for all $i\in\mathcal{V}_M$, $M_i=m>0$, $D_i=d>0$, $\alpha_i=1$. For DPIAC, assume 
that the topology of the communication network is the same as the one of the power system such that $l_{ij}=K_{ij}$ for all $(i,j)\in\mathcal{E}$. 
\end{assumption}
\par The frequency dependent nodes are excluded from the model by this assumption while the nodes 
of the passive power loads can be involved into the model in this assumption by Kron Reduction \cite{P.M.Anderson}. From the practical point of view, the analysis with Assumption  \ref{assumption_a1} is valuable because it provides us the insight on how to improve the transient behavior by tuning the control coefficients. For the general case without the restriction of this assumption, in which the model includes the frequency dependent nodes, we resort to simulations in Section \ref{section:simulationstudy}.
\par 
Since $\theta_{ij}=\theta_i-\theta_j$ is usually small for relatively small $P_i$ compared to the line capacity in practice, we approximate $\sin{\theta_{ij}}$ by $\theta_{ij}$ as in e.g., \cite{Zhao2015,EAGC} to focus on the transient performance. With Assumption \ref{assumption_a1}, rewriting (\ref{chapterModel:eq:system1}) into a vector form by replacing $\sin{\theta_{ij}}$ by $\theta_{ij}$, we obtain 
\begin{subequations}\label{linearized:system}
 \begin{align}
 \dot{\bm\theta}&=\bm\omega\\
  \bm M\dot{\bm \omega}&=-\bm L\bm\theta-\bm D\bm\omega+\bm B\bm w+\bm u,
 \end{align}
\end{subequations}
where $\bm \theta=\text{col}(\theta_i)\in\mathbb{R}^{n}$, $n$ denotes the number of nodes in the network,  $\bm\omega=\text{col}(\omega_i)\in\mathbb{R}^n$, $\bm M=\text{diag}(M_i)\in\mathbb{R}^{n\times n}$, $\bm L\in\mathbb{R}^{n\times n}$ is the Laplacian 
 matrix of the network, 
 $\bm D=\text{diag}(D_i)\in\mathbb{R}^{n\times n}$, $\bm u=\text{col}(u_i)\in\mathbb{R}^{n\times n}$. The disturbances of $\bm P=\text{col}(P_i)\in\mathbb{R}^n$ have been modeled by $\bm B \bm w$ as the inputs with $\bm B\in\mathbb{R}^{n\times n}$ and $\bm w\in\mathbb{R}^n$ as in Definition \ref{definition_H2norm}. Here, $\text{col}(\cdot)$ denotes 
the column vector of the indicated elements and $\text{diag}{(\beta_i)}$ denotes 
a diagonal matrix $\bm\beta=\text{diag}{(\{\beta_i,i\cdots n\})}\in\mathbb{R}^{n\times n}$ with $\beta_i\in\mathbb{R}$. 
Denote the identity matrix by $\bm I_n\in\mathbb{R}^{n\times n}$ and the $n$ dimensional vector with all elements equal to one by $\bm 1_n$.

\par 
The transient performance of $\bm \omega(t)$ and $\bm u(t)$ are measured by the $\mathcal{H}_2$ norm of the corresponding transfer functions with input $\bm w$ and output $\bm y=\bm \omega$ and $\bm y=\bm u$ respectively. 
The squared $\mathcal{H}_2$ norms are denoted by $||\bm G_i(\bm \omega,\bm w)||^2_2$ and $||\bm G_i(\bm u,\bm w)||^2_2$ where the sub-index $i=c~\text{or}~d$ which refers to the centralized method GBPIAC or the distributed method DPIAC.



\section{The transient performance analysis}\label{chapter5:Sec:coherence}

In this section, we calculate the $\mathcal{H}_2$ norms of the frequency deviation and of the control cost for GBPIAC and DPIAC. 
\par 
\begin{lemma}\label{lemma_performance}
 For a symmetric Laplacian matrix $\bm L\in\mathbb{R}^{n\times n}$,
 there exist an invertible matrix $\bm Q\in\mathbb{R}^{n\times n}$ such that 
\begin{subequations}
 \begin{align}
  \bm Q^{-1}=\bm Q^T,\\
  \bm Q^{-1}\bm L\bm Q=\bm \Lambda,\\
  \bm Q_1=\frac{1}{\sqrt{n}}\bm 1_n,
 \end{align}
\end{subequations}
where $\bm Q=[\bm Q_1,\cdots,\bm Q_n]$, $\bm\Lambda=\text{diag}(\lambda_i)\in\mathbb{R}^{n\times n}$, $\bm Q_i\in\mathbb{R}^n$ is the normalized eigenvector of $\bm L$ corresponding to eigenvalue $\lambda_i$, thus 
$\bm Q^T_i\bm Q_j=0$ for $i\neq j$. 
Because $\bm L\bm 1_n=\bm 0$, $\lambda_1= 0$ is one of the eigenvalues with normalized eigenvector $\bm Q_1$. 
\end{lemma}
\par 
We study the 
transient performance of $\bm \omega$ and $\bm u$ of GBPIAC and DPIAC in subsection \ref{chapter5:subsec:performance_GBPIAC} and \ref{chapter5:subsec:performance_analysis_DPIAC} respectively by 
calculating the corresponding $\mathcal{H}_2$ norm. In addition, 
for DPIAC, we also calculate a $\mathcal{H}_2$ norm 
which measures the coherence of the marginal costs. The performance 
of GBPIAC and DPIAC will be compared in subsection \ref{chapter5:subsec:comparison}.

\subsection{Transient performance analysis for GBPIAC}\label{chapter5:subsec:performance_GBPIAC}
By Assumption \ref{assumption_a1}, we derive the control input $u_i=\frac{1}{n}k_1\xi_s$ at node $i$ as  in (\ref{eq:GBPIAC}). 
With the notations of section \ref{chapter5:Sec:problem} and Assumption \ref{assumption_a1}, and $k_2=4k_1$,
we obtain from (\ref{linearized:system}) and (\ref{eq:GBPIAC}) the 
closed-loop system of GBPIAC in a vector form as follows.
\begin{subequations}\label{chapter5:eq:closed_GBPIAC}
 \begin{align}
  \dot{\bm\theta}&=\bm\omega,\label{chapter5:closed_GBPIACa}\\
  m\bm I_n\dot{\omega}&=-\bm L\bm\theta-d\bm I_n\omega+\frac{4k_1\xi_s}{n}\bm 1_n+\bm B\bm w,\label{chapter5:closed_GBPIACb}\\
  \dot{\eta}_s&=d\bm 1_n^{T}\bm\omega,\\
  \dot{\xi}_s&=-k_1m\bm 1_n^T\bm\omega-k_1\eta_s-4k_1\xi_s, 
 \end{align}
\end{subequations}
where $\eta_s\in\mathbb{R}$ and $\xi_s\in\mathbb{R}$.

For the transient performance of $\bm\omega(t), \bm u(t)$ in GBPIAC, the following theorem can be proved.
\begin{theorem}\label{theorem1_h2norm}
Consider the closed-loop system (\ref{chapter5:eq:closed_GBPIAC}) of GBPIAC with $\bm B=\bm I_n$. The squared $\mathcal{H}_2$ norm of the frequency deviation $\bm\omega$ and of the control inputs $\bm u$ are,
 \begin{subequations}
  \begin{align}
 ||\bm G_c(\bm\omega,\bm w)||^2_2&=\frac{n-1}{2md}+\frac{d+5mk_1}{2m(2k_1m+d)^2},\label{chapter5:expected:frequency_GBPIAC}\\
 ||\bm G_c(\bm u,\bm w)||^2_2&=\frac{k_1}{2}.\label{chapter5:expect_u_GBPIAC}
  \end{align}
 \end{subequations}  
\end{theorem}
\emph{Proof:} With the linear transform $\bm x_1=\bm Q^{-1}\bm \theta, \bm x_2=\bm Q^{-1}\bm \omega$ where $\bm Q$ is defined 
in Lemma \ref{lemma_performance}, we derive from (\ref{chapter5:eq:closed_GBPIAC}) that
\begin{subequations}
 \begin{align*}
\dot{\bm x}_1&=\bm x_2,\\
\dot{\bm x}_2&=-\frac{1}{m}\bm\Lambda \bm x_1-\frac{d}{m}\bm I_n \bm x_2+\frac{4k_1\xi_s}{mn}\bm Q^{-1}1_n+\frac{1}{m}\bm Q^{-1}\bm w,\\
\dot{\eta}_s&=d\bm 1_n^T\bm Q\bm x_2,\\
\dot{\xi}_s&=-k_1m\bm 1_n^T\bm Q\bm x_2-k_1\eta_s-4k_1\xi_s,
 \end{align*}
\end{subequations}
where $\bm\Lambda$ is the diagonal matrix defined in Lemma \ref{lemma_performance}. 
Since $\bm 1_n$ is an eigenvector of $\bm L$ corresponding to $\lambda_1=0$, we obtain $\bm Q^{-1}\bm 1_n=[\sqrt{n},0,\cdots,0]^T$.
Thus the components of $\bm x_1$ and $\bm x_2$ 
can be decoupled as  
\begin{subequations}\label{chapter5:decoupled1}
 \begin{align}
  \dot{x}_{11}&=x_{21},\\
  \dot{x}_{21}&=-\frac{d}{m}x_{21}+\frac{4k_1}{m\sqrt{n}}\xi_s+\frac{1}{m}\bm Q_{1}^T\bm w,\\
  \dot{\eta}_s&=d\sqrt{n}x_{21},\\
  \dot{\xi}_s&=-k_1m\sqrt{n}x_{21}-k_1\eta_s-4k_1\xi_s
 \end{align}
\end{subequations}
and for $i=2,\cdots,n$, 
\begin{subequations}\label{chapter5:decoupled2}
 \begin{align}
  \dot{x}_{1i}&=x_{2i},\\
  \dot{x}_{2i}&=-\frac{\lambda_i}{m}x_{1i}-\frac{d}{m}x_{2i}+\frac{1}{m}\bm Q_{i}^T\bm w,
 \end{align}
\end{subequations}
We rewrite the decoupled systems of (\ref{chapter5:decoupled1}) and (\ref{chapter5:decoupled2}) 
in the general form as (\ref{Appendix:generalform}) with 
\[\bm x=\begin{bmatrix}
     \bm x_1\\
     \bm x_2\\
     \eta_s\\
     \xi_s
    \end{bmatrix},
 \bm A=\begin{bmatrix}
    \bm 0&\bm I_n&\bm 0&\bm 0\\
    -\frac{\bm\Lambda}{m}&-\frac{d}{m}\bm I_n&\bm 0&\frac{4k_1}{m\sqrt{n}}\bm v\\
    \bm 0&d\sqrt{n}\bm v^T&\bm 0&\bm 0\\
    \bm 0&-k_1m\sqrt{n}\bm v^T&-k_1&-4k_1
   \end{bmatrix},
 \bm \tilde{B}=\begin{bmatrix}
    \bm 0\\
    \frac{\bm Q^{-1}}{m}\\
    \bm 0\\
    \bm 0
   \end{bmatrix},
\]
where $\bm v^T=[1,0,\cdots,0]\in\mathbb{R}^n$. The $\mathcal{H}_2$ norm of a state variable e.g., the 
frequency deviation and the control cost, can be determined by setting the output $y$ as that 
state variable. 
Because the closed-loop 
system (\ref{chapter5:eq:closed_GBPIAC}) is asymptotically stable, $\bm A$ is Hurwitz regardless the rotations of the phase angle $\bm\theta$.

For the transient performance of $\bm \omega(t)$, setting $\bm y=\bm\omega=\bm Q\bm x_2$ and $\bm C=[\bm 0,\bm Q,\bm 0,\bm 0]$, 
we obtain the observability Grammian $\bm Q_o$ of $(\bm C,\bm A)$ (\ref{Lyapunov:equation}) in
the form,
\[\bm Q_o=\begin{bmatrix}
       \bm Q_{o11}&\bm Q_{o12}&\bm Q_{o13}&\bm Q_{o14}\\
       \bm Q_{o12}^T&\bm Q_{o22}&\bm Q_{o23}&\bm Q_{o24}\\
       \bm Q_{o13}^T&\bm Q_{o23}^T&\bm Q_{o33}&\bm Q_{o34}\\
       \bm Q_{o14}^T&\bm Q_{o24}^T&\bm Q_{o34}^T& Q_{o44}\\
      \end{bmatrix}.
\]
Thus, 
\begin{eqnarray}
 \|\bm G_c(\bm \omega,\bm w)\|^2_2=\text{tr}(\tilde{\bm B}^T\bm Q_o\tilde{\bm B})=\frac{\text{tr}(\bm Q\bm Q_{o22}\bm Q^T)}{m^2}=\frac{\text{tr}(\bm Q_{o22})}{m^2}. \label{trace:equivalent}
\end{eqnarray}
Because
\[\bm C^T\bm C=\begin{bmatrix}
        \bm 0& \bm 0 &\bm 0 &\bm 0\\
        \bm 0&\bm I_n&\bm 0 &\bm 0\\
        \bm 0& \bm 0 &\bm 0 &\bm 0\\
        \bm 0& \bm 0 &\bm 0 &\bm 0\\
       \end{bmatrix}, 
\]
the diagonal elements $\bm Q_{o22}(i,i)$ of $\bm Q_{o22}$ can be calculated by solving the observability Gramian $\tilde{\bm Q}_i$ of $(\bm C_i,\bm A_i)$
which satisfies
\begin{eqnarray*}
 \tilde{\bm Q}_1\bm A_1+\bm A^T_1\tilde{\bm Q}_1+\bm C_1^T\bm C_1=0
\end{eqnarray*}
where
\[\bm A_1=\begin{bmatrix}
     0 &1 &0 &0\\
     0 &-\frac{d}{m}&0&\frac{4k_1}{m\sqrt{n}}\\
     0 &d\sqrt{n}  &0 &0 \\
     0 &-k_1m\sqrt{n}&-k_1&-4k_1
    \end{bmatrix},
  \bm C_1^T=\begin{bmatrix}
         0\\
         1\\
         0\\
         0
        \end{bmatrix}
\]
and 
\begin{eqnarray*}
 \tilde{\bm Q}_i\bm A_i+\bm A^T_i\tilde{\bm Q}_i+\bm C_i^T\bm C_i=\bm 0,~i=2,\cdots,n,
\end{eqnarray*}
where
\[\bm A_i=\begin{bmatrix}
       0&1\\
       -\frac{\lambda_i}{m}&-\frac{d}{m}
      \end{bmatrix}, 
  \bm C_i^T=\begin{bmatrix}
         0\\
         1\\
        \end{bmatrix}.
\]
In this case, the diagonal elements of $\bm Q_{o22}$ satisfy $\bm Q_{o22}(i,i)=\tilde{\bm Q_i}(2,2)$ for $i=1,\cdots,n$. We thus derive
from the observability Gramian $\tilde{\bm Q}_i$ that
\begin{eqnarray}
 \text{tr}(\tilde{\bm B}^T\bm Q_o\tilde{\bm B})=\frac{n-1}{2md}+\frac{d+5mk_1}{2m(2k_1m+d)^2}, \label{chapter5:procedure}
\end{eqnarray}
which yields (\ref{chapter5:expected:frequency_GBPIAC}) directly. 
Similarly by setting $y=u_s(t)=4k_1\xi_s(t)$ and $\bm C=[\bm 0,\bm 0,0,4k_1]$, we derive the norm of $u_s(t)$ as 
\begin{eqnarray}
 \|\bm G_c(u_s,\bm w)\|^2_2=\frac{k_1n}{2}.\label{expected_Us}
\end{eqnarray}
With $u_i=\frac{u_s}{n}$ for $i=1,\cdots,n$ and $\|\bm u(t)\|^2=\bm u(t)^T\bm u(t)=\sum_i^n{u_i^2(t)}$, 
we further derive(\ref{chapter5:expect_u_GBPIAC}) for the control cost.
\hfill $\Box$
\par 

The norm of $\bm \omega(t)$ in (\ref{chapter5:expected:frequency_GBPIAC}) includes two terms. 
The first one describes \emph{the relative deviations} which depend on the primary control, and the second one describes \emph{the overall frequency deviation} of which the suppression is the task of the secondary control. From the proof of Theorem \ref{theorem1_h2norm}, 
it can be observed that the relative frequency deviations are derived
from the eigen-direction of nonzero eigenvalues, and the
overall frequency deviation from the eigen-direction of the zero eigenvalue.
\par
\begin{remark}
It is demonstrated by Theorem \ref{theorem1_h2norm} that the topology of the 
network has no influence on the norm of $\bm\omega(t)$ and $\bm u(t)$. This is because of Assumption \ref{assumption_a1} and the identical strength of the disturbances at all the nodes with $\bm B=\bm I_n$.
\end{remark}
\begin{remark}\label{remark2GB}

The overall frequency deviation depends on $k_1$ while the relative frequency deviation is independent of $k_1$. Hence, the frequency deviation cannot be 
suppressed to an arbitrary small positive value. When the overall frequency deviation caused by the power imbalance dominates the relative frequency deviation, a large $k_1$ can accelerate the restoration of the frequency. This will be 
further described in Section \ref{section:simulationstudy}. However, a large $k_1$ leads to a high control cost. Hence there is a trade-off between the overall frequency deviation suppression and control cost, which is determined by $k_1$. 
\end{remark}
\par 
%
%
%

Theorem \ref{theorem1_h2norm} with $\bm B=\bm I_n$ includes the assumption that all the disturbances are independent and of identical strength. The following 
theorem describes the impact of the control coefficients with a general $\bm B$ 
where the disturbances are correlated with non-identical strength.
\begin{theorem}\label{theorem1}
Consider the closed-loop system (\ref{chapter5:eq:closed_GBPIAC}) of GBPIAC 
with a positive definite $\bm B\in\mathbb{R}^{n\times n}$, the squared $\mathcal{H}_2$ norm of the frequency deviation $\bm \omega$ and of the control inputs $\bm u$ satisfy
\begin{subequations}\label{inequalities1}
\begin{align}
\gamma_{\min}^2G_c&\leq ||\bm G_c(\bm\omega,\bm w)||^2_2\leq \gamma_{\max}^2G_c,\label{inequalities01}\\
\frac{k_1}{2}\gamma_{\min}^2&\leq ||\bm G_c(\bm u,\bm w)||^2_2\leq \frac{k_1}{2}\gamma_{\max}^2,\label{inequalities02}
\end{align}
\end{subequations}
where  $\gamma_{\min}$ and $\gamma_{\max}$ are the smallest and the largest 
eigenvalue of $\bm B$, and 
\begin{align*}
G_c=\big(\frac{n-1}{2md}+\frac{d+5mk_1}{2m(2k_1m+d)^2}\big).
\end{align*}
\end{theorem}
\emph{Proof:} 
Since $\bm B>0$, then there exist $\bm U$ and $\bm \Gamma$ such that $\bm B=\bm U^T\bm \Gamma\bm U$ where $\bm U$ is an orthogonal matrix and $\bm \Gamma$ is a diagonal matrix with the diagonal elements being the eigenvalues of $\bm B$, which are all strictly positive. With the decomposition of $\bm B$,  from (\ref{trace:equivalent}) we derive
\begin{subequations}
\begin{align*}
\text{trac}\big(\bm B^T\bm Q\bm Q_{o22}Q\bm B\big)&=\text{trac}\big(\bm U^T\bm \Gamma\bm U\bm Q\bm Q_{o22}\bm Q^T\bm U^T\bm \Gamma\bm U\big)\\
&=\text{trac}\big(\bm \Gamma\bm U\bm Q\bm Q_{o22}\bm Q^T\bm U^T\bm \Gamma\big)\\
&\leq \gamma_{\max}^2\text{trac}\big(\bm U\bm Q\bm Q_{o22}\bm Q^T\bm U^T\big)\\
&=\gamma_{\max}^2\text{trac}\big(\bm Q_{o22}\big),
\end{align*}
\end{subequations}
Similarly, we derive 
\begin{subequations}
\begin{align*}
\gamma_{\min}^2\text{trac}\big(\bm Q_{o22}\big)\leq \text{trac}\big(\bm B^T\bm Q\bm Q_{o22}\bm Q^T\bm B \big).
\end{align*}
\end{subequations}
From the above inequalities and (\ref{trace:equivalent}), we can easily follow
the procedure as in the proof of Theorem \ref{theorem1_h2norm} to obtain the trace of $Q_{o22}$ and further derive
the inequalities in (\ref{inequalities01}). A similar procedure is conducted to obtain 
the inequalities in (\ref{inequalities02}). \hfill $\Box$
\par 
~~
\par 
It is demonstrated by Theorem \ref{theorem1} that when the disturbances are correlated with non-identical
strength, the impact of $k_1$ on the norm is similar as in the case with identical strength of disturbances.
\par

\subsection{Transient performance analysis for DPIAC}\label{chapter5:subsec:performance_analysis_DPIAC}

With Assumption \ref{assumption_a1} and $k_2=4k_1$, we derive 
the closed-loop system of DPIAC from (\ref{linearized:system}) and (\ref{eq:DPIAC}) as 
\begin{subequations}\label{chapter5:system2}
 \begin{align}
 \dot{\bm \theta}&=\bm\omega,\\
 m\bm I_n\dot{\bm\omega}&=-\bm L\bm\theta-d\bm I_n\bm\omega+4k_1\bm\xi+\bm B\bm w\\
 \dot{\bm\eta}&=\bm D\bm\omega+4k_1k_3\bm L\bm \xi,\\
 \dot{\bm\xi}&=-k_1\bm M\bm\omega-k_1\bm\eta-4k_1\bm\xi,
 \end{align}
\end{subequations}
where $\bm\eta=\text{col}(\eta_i)\in\mathbb{R}^n$ and $\bm\xi=\text{col}(\xi_i)\in\mathbb{R}^n$. Note that $\bm L=(L_{ij})\in \mathbb{R}^{n\times n}$ is the
weighted Laplacian matrix of the power network and also of the communication network.
Because the differences of the marginal costs can be fully represented by $4k_1\bm L\bm\xi(t)$, we use the squared norm of $(4k_1\bm L\bm\xi(t))$ to measure the coherence of the marginal costs in DPIAC. 
Denote the squared $\mathcal{H}_2$ norm of the transfer matrix of (\ref{chapter5:system2}) with output $\bm y=4k_1\bm L\bm\xi$ by 
$||G_d(4k_1L\xi,\bm w)||^2$.
In this subsection, we also calculate the squared $\mathcal{H}_2$ norm of $4k_1\bm L\bm\xi$ as an additional metric for the influence
of $k_3$ on the control cost.

The following theorem states the $\mathcal{H}_2$ norms of the frequency deviation,  the control cost and 
the coherence of the marginal costs in DPIAC.
\begin{theorem}\label{theorem2}
 Consider the closed-loop system (\ref{chapter5:system2}) of DPIAC with $\bm B=\bm I_n$, the squared $\mathcal{H}_2$ norm of $\bm\omega(t)$, $\bm u(t)$ and $4k_1\bm L\bm\xi$ are 
 \begin{subequations}
  \begin{align}
   \|\bm G_d(\bm\omega,\bm w)\|^2_2&=\frac{1}{2m}\sum_{i=2}^n\frac{b_{1i}}{e_{i}}+\frac{d+5mk_1}{2m(2k_1m+d)^2},\label{special_omega}\\
   \|\bm G_d(\bm u,\bm w)\|^2_2&=\frac{k_1}{2}+\sum_{i=2}^{n}\frac{b_{2i}}{e_{i}},\label{special_u}\\
   \|\bm G_d(4k_1\bm L\bm\xi,\bm w)\|^2_2&=\sum_{i=2}^{n}\frac{\lambda_i^2b_{2i}}{m^2e_{i}}, \label{special_marginal}
  \end{align}
 \end{subequations}
 where 
 \begin{subequations}
  \begin{align*}
    b_{1i}&=\lambda_i^2(4k_1^2k_3m-1)^2+4dmk_1^3\\
  &+k_1(d+4k_1m)(4d\lambda_ik_1k_3+5\lambda_i+4dk_1),\\
    b_{2i}&=2dk_1^3(d+2k_1m)^2+2\lambda_ik_1^4m^2(4k_1k_3d+4),\\
  e_{i}&=d\lambda_i^2(4k_1^2k_3m-1)^2+16d\lambda_ik_1^4k_3m^2+d^2\lambda_ik_1\\
  &+4k_1(d+2k_1m)^2(dk_1+\lambda_i+d\lambda_i k_1k_3).
  \end{align*}
 \end{subequations}
\end{theorem}
\emph{Proof:} Let $\bm Q\in\mathbb{R}^{n\times n}$ be defined as in Lemma \ref{lemma_performance} and let $\bm x_1=\bm Q^{-1}\bm\theta, \bm x_2=\bm Q^{-1}\bm\omega, \bm x_3=\bm Q^{-1}\bm\eta, \bm x_4=\bm Q^{-1}\bm\xi$, 
we obtain the closed-loop system in the general form as (\ref{Appendix:generalform})
with 
{\small
\[\bm x=\begin{bmatrix}
     \bm x_1\\
     \bm x_2\\
     \bm x_3\\
     \bm x_4
    \end{bmatrix},
 \bm A=\begin{bmatrix}
    \bm 0&\bm I_n&\bm 0&\bm 0\\
    -\frac{1}{m}\bm\Lambda&-\frac{d}{m}\bm I_n&\bm 0&\frac{4k_1}{m}\bm I_n\\
    \bm 0&d\bm I_n&\bm 0&4k_1k_3\bm\Lambda\\
    \bm 0&-k_1m\bm I_n&-k_1\bm I_n&-4k_1\bm I_n
   \end{bmatrix},
 \bm B=\begin{bmatrix}
    \bm 0\\
    \frac{\bm Q^{-1}}{m}\\
    \bm 0\\
    \bm 0
   \end{bmatrix},
\]}
where $\bm\Lambda$ is the diagonal matrix defined in Lemma \ref{lemma_performance}. 
Each of the block
matrices in the matrix $\bm A$ is either the zero matrix or a diagonal matrix, so the components of the vector $\bm x_1\in\mathbb{R}^n,~\bm x_2\in\mathbb{R}^n,~\bm x_3\in\mathbb{R}^n,~\bm x_4\in\mathbb{R}^n$ 
can be decoupled. 

With the same method for obtaining (\ref{chapter5:procedure}) in the proof of Theorem \ref{theorem1_h2norm}, setting  $\bm y=\bm\omega=Qx_2$, $\bm C=[\bm 0,\bm Q,\bm 0,\bm 0]$, we derive (\ref{special_omega}) for $\bm \omega(t)$. Then, setting $\bm y(t)=\bm u(t)=4k_1\bm\xi(t)$ and $\bm C=[\bm 0,\bm 0,\bm 0,4k_1\bm Q]$,
we derive (\ref{special_u}) for the norm of $\bm u(t)$. Finally for the coherence measurement of the marginal cost, setting $\bm y=4k_1\bm L\bm\xi$ and 
$\bm C=[\bm 0,\bm 0,\bm 0,4k_1\bm L\bm Q]$, we derive (\ref{special_marginal}). 
\hfill $\Box$
\par 

Similar as Theorem \ref{theorem1}, for a positive definite $\bm B\in\mathbb{R}^{n\times n}$, we obtain the following theorem.
\begin{theorem}
Consider the closed-loop system (\ref{chapter5:system2}) of DPIAC with a positive definite $\bm B\in\mathbb{R}^{n\times n}$, the squared $\mathcal{H}_2$ norm of $\bm\omega(t)$, $\bm u(t)$ and $4k_1\bm L\bm\xi$ satisfy 
\begin{subequations}
\begin{align*}
\hspace{-20pt}
\gamma_{\min}^2G_d&\leq\|\bm G_d(\bm\omega,\bm w)\|^2_2\leq \gamma_{\max}^2G_d,\\
\hspace{-20pt}
\gamma_{\min}^2(\frac{k_1}{2}+\sum_{i=2}^{n}\frac{b_{2i}}{e_{i}})&\leq \|\bm G_d(\bm u,\bm w)\|^2_2\leq \gamma_{\max}^2(\frac{k_1}{2}+\sum_{i=2}^{n}\frac{b_{2i}}{e_{i}})\\
\hspace{-20pt}
\gamma_{\min}^2\sum_{i=2}^{n}\frac{\lambda_i^2b_{2i}}{m^2e_{i}}&\leq \|\bm G_d(4k_1\bm L\bm\xi,\bm w)\|^2_2\leq\gamma_{\max}^2\sum_{i=2}^{n}\frac{\lambda_i^2b_{2i}}{m^2e_{i}}
\end{align*}
\end{subequations}
where  $\gamma_{\min}$ and $\gamma_{\max}$ are the smallest and the largest 
eigenvalue of $\bm B$, $b_{1i},b_{2i}$ and $e_i$ are defined in Theorem \ref{theorem2} and 
\begin{align*}
G_d=\frac{1}{2m}\sum_{i=2}^n\frac{b_{1i}}{e_{i}}+\frac{d+5mk_1}{2m(2k_1m+d)^2}. 
\end{align*}
\end{theorem}
The proof of this theorem follows that of Theorem \ref{theorem1}. 
Hence, similar as in GBPIAC, when the disturbances are correlated with non-identical strength, the impact of $k_1$ and $k_3$ on the norms are similar as in the case with identical strength of disturbances.
\par

Based on Theorem \ref{theorem2}, we analyze the impact of $k_1$ and $k_3$ on the norms by focusing on 1) the frequency deviation, 2) control cost and 3) coherence of the marginal costs.

\subsubsection{The frequency deviation}\label{chapter5:subsec:frequency}
We first pay attention to 
the influence of $k_1$ when $k_3$ is fixed. The norm of $\bm \omega$ also includes two terms in (\ref{special_omega}) where the first one describes the 
relative frequency oscillation and the second one describes the overall frequency 
deviation. The overall frequency deviation decreases 
inversely as $k_1$ increases. Hence, when the overall frequency deviation dominates the relative frequency deviation, the convergence can also be accelerated by a large $k_1$ as analyzed in 
Remark \ref{remark2GB} for GBPIAC. From (\ref{special_omega}), we derive 
\begin{eqnarray}\label{k1infinity}
\lim_{k_1\rightarrow\infty}{||\bm G_d(\bm\omega,\bm w)||^2_2}=\frac{1}{2m}\sum_{i=2}^{n}\frac{\lambda_i^2k_3^2}{d\lambda_i^2k_3^2+d(1+2\lambda_ik_3)},
\end{eqnarray}
which indicates that even with a large $k_1$, the frequency deviations cannot be decreased anymore when $k_3$ is nonzero. \emph{So similar to GBPIAC, the relative frequency deviation cannot be suppressed to an arbitrary small positive value in DPIAC}. However, when $k_3=0$, DPIAC reduces to DecPIAC (\ref{decentralized}), thus
\begin{eqnarray}
\|\bm G_d(\bm\omega,\bm w)\|^2_2\sim O(k_1^{-1}). \label{chapter5:omega_norm2}
\end{eqnarray} 
\begin{remark}\label{chapter5:remark2}
By DecPIAC, it follows from (\ref{chapter5:omega_norm2}) that if all the nodes are equipped with the secondary frequency controllers, the frequency deviation can be 
controlled to any prespecified range. However, it results in a high control cost for the entire network. In addition, the configuration of $k_1$ 
is limited by the response time of the actuators. 
\end{remark}

\begin{remark}
 This analysis is based on Assumption \ref{assumption_a1} which requires that each node in the network is equipped with a secondary frequency controller. 
 However, for the power systems without all the nodes equipped with the controllers, the disturbance 
 from the node without a controller must be compensated by the other nodes with controllers. In that case,
 the equilibrium of the system is changed and the oscillation can never be avoided even when the controllers are sufficiently sensitive to the disturbances. 
\end{remark}
 
When $k_1$ is fixed, it can be 
easily observed from (\ref{special_omega}) that the order of $k_3$ in the term $b_{1i}$ is 2 which is the same as in the term $e_i$, thus $k_3$ 
has little influence on the frequency deviation.

\subsubsection{The control cost}
We first analyze the influence of $k_1$ on the cost and then the influence of $k_3$. For any $k_3\geq 0$, 
we derive from (\ref{special_u}) that 
 \begin{equation*}
 \|\bm G_d(\bm u,\bm w)\|^2_2\sim O(k_1),
 \end{equation*}
which indicates that the control cost increases as $k_1$ increases. Recalling 
the impact of $k_1$ on the overall frequency deviation in (\ref{special_omega}),
we conclude that minimizing the control cost always conflicts with minimizing the frequency deviation. \emph{Hence, a trade-off 
should be determined to obtain the desired frequency deviation with an acceptable control cost. }

Next, we analyze how $k_3$ influences the control cost. From (\ref{special_u}), we 
obtain that 
\begin{eqnarray}\label{k3infinity}
{\|\bm G_d(\bm u,\bm w)\|^2_2}\sim\frac{k_1}{2}+O(k_1k_3^{-1}),
\end{eqnarray}
where the second term is positive. It shows that
the control cost decreases as $k_3$ increases due to the accelerated consensus speed of the marginal costs. This will be 
further discussed in the next subsubsection on the coherence of the marginal 
costs. 
Note that $k_3$ has little influence on the frequency deviation. 
\emph{Hence the control cost can be decreased by $k_3$ without increasing the frequency deviation much. }

\subsubsection{The coherence of the marginal costs in DPIAC}\label{chapter5:subsubsec:coherence}

We measure the coherence of the marginal costs by the norm of 
$\|\bm G_d(4k_1\bm L\bm\xi,\bm w)||^2$. 
From (\ref{special_marginal}), we obtain 
\begin{eqnarray*}
\|\bm G_d(4k_1\bm L\bm\xi,\bm w)\|^2_2=O(k_3^{-1}),
\end{eqnarray*}
which indicates that the  difference of the marginal costs decreases as $k_3$ increases. \emph{Hence, this analytically 
confirms that the consensus speed can be 
increased by increasing $k_3$.} 

\begin{remark}
In practice, similar to $k_1$, the configuration of $k_3$ depends on the communication devices and cannot be arbitrarily large. In addition, the communication delay and noise also influence the transient performance, which still needs further investigation. 
\end{remark}

%

\subsection{Comparison of the GBPIAC and DPIAC control laws}\label{chapter5:subsec:comparison}

With a positive $k_1$, we can easily obtain from (\ref{chapter5:expect_u_GBPIAC}, \ref{special_u}) that 
\begin{eqnarray}
\|\bm G_c(\bm u,\bm w)\|<\|\bm G_d(\bm u,\bm w)||,
\end{eqnarray}
which is due to the differences of the marginal costs. The difference in the control cost between 
these two control laws can be decreased by accelerating the consensus of the marginal costs as explained in the previous subsection. 
From (\ref{special_omega}) and (\ref{special_u}) we derive that
\begin{subequations}
 \begin{align*}
\hspace{-18pt}
\lim_{k_3\rightarrow\infty}{\|\bm G_d(\bm\omega,\bm w)\|^2_2}&=\frac{n-1}{2md}+\frac{d+5mk_1}{2m(2k_1m+d)^2}=||\bm G_c(\bm\omega,\bm w)||^2,\\
\hspace{-18pt}
\lim_{k_3\rightarrow\infty}{\|\bm G_d(\bm u,\bm w)\|^2_2}&=\frac{k_1}{2}=\|\bm G_c(\bm u,\bm w)\|^2. 
 \end{align*}
\end{subequations}
Hence, as $k_3$ goes to infinity, the transient performance of DPIAC converges to that of GBPIAC.

\section{Simulations}\label{section:simulationstudy}
\begin{figure}[ht]
 \begin{center}
\includegraphics[height=185pt,width=250pt]{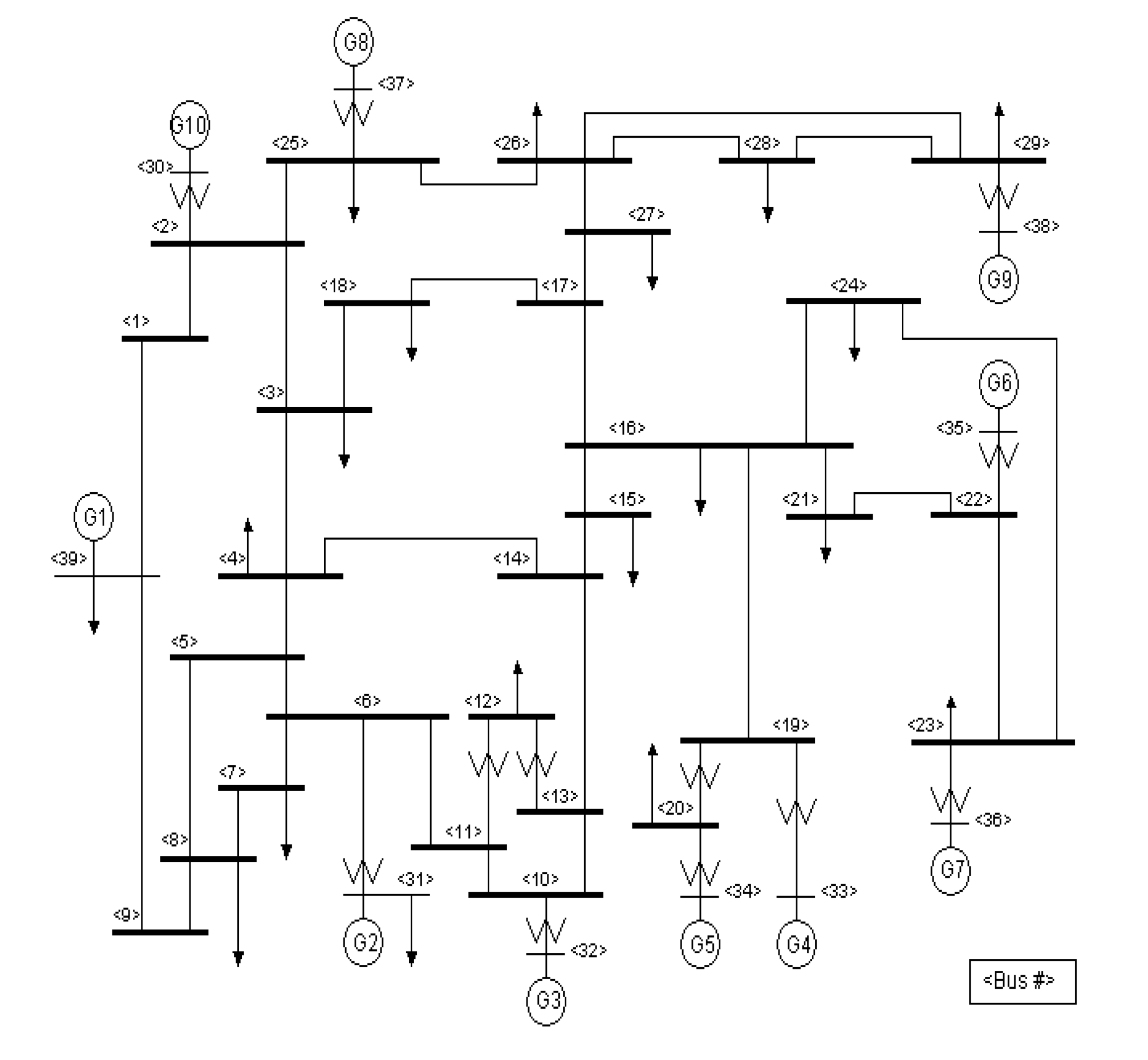}
 \caption{IEEE 39-bus test power system}
  \label{fig.IEEE39_graph}
 \end{center}
\end{figure}

In this section, we numerically verify the analysis of the transient performance of DPIAC using the IEEE 39-bus system as shown in Fig.~\ref{fig.IEEE39_graph} with the Power System Analysis Toolbox (PSAT)~\cite{Milano2008}. We compare the performance of DPIAC with that of GBPIAC. For a comparison of DPIAC with the traditional 
control laws, see \cite{PIAC3,MLPIAC}. The system consists
of 10 generators, 39 buses, which serves a total load of
about 6 GW. As in \cite{MLPIAC}, we change the buses which are neither connected to synchronous machines nor to power loads into frequency dependent buses. Hence $\mathcal{V}_{M}=\{G1,G2,G3,G4,G5,G6,G7,G8,G9,G10\}$, $\mathcal{V}_{P}=\{30,31,32,33,34,35,36,37,38,39\}$ 
and the other nodes are in set $\mathcal{V}_{F}$. The nodes 
in $\mathcal{V}_M\cup\mathcal{V}_F$ are all equipped with secondary frequency controllers such that $\mathcal{V}_K=\mathcal{V}_M\cup\mathcal{V}_F$. Because the voltages are constants, the angles of the synchronous machine and the bus have the same dynamics \cite{Ilic2000}.  
Except the control gain coefficients $k_1$ and $k_3$, all the parameters of the power system, including the control prices, damping coefficients and constant voltages are identical to those in the simulations in \cite{MLPIAC}. The communication topology are the same as the one of the power network and we set 
 $l_{ij}=1$ for the communication if node $i$ and $j$ are connected. We remark that with
 these configurations of the parameters, Assumption \ref{assumption_a1} is not satisfied in the simulations. 
We first verify the impact of $k_1$ and $k_3$ on the transient performance in the deterministic system where the disturbance is modeled by a step-wise increase of load, then in a stochastic system with the interpretation of the $\mathcal{H}_2$ norm as the limit of the variance of the output. 
\subsection{In the deterministic system}
We analyze the impact of the control gain coefficients
on the performance on the deterministic system where the disturbances are 
step-wise increased power loads by 66 MW at nodes 4, 12 and 20 at time $t=5$ seconds. This step-wise disturbance lead the overall frequency  to dominate the relative frequencies as described in Remark \ref{remark2GB}, which illustrates
the function of the secondary frequency control. The system behavior following the disturbance also show us how the convergence
of the state can be accelerated by tuning $k_1$ and $k_3$ monotonically.
We calculate the following two metrics
\begin{align*}
S=\int_0^{T_0}\bm\omega^T(t)\bm\omega(t)dt,~\text{and}  
~C=\frac{1}{2}\int_0^{T_0}\bm u^T(t)\bm\alpha\bm u(t)dt,
\end{align*}
to measure the performance of $\bm \omega(t)$ and $\bm u(t)$ during the transient phase,
where $T_0=40$, $\bm\omega=\text{col}(\omega_i)$ for $i\in\mathcal{V}_M\cup\mathcal{V}_F$ , $\bm u=\text{col}(u_i)$ for $i\in\mathcal{V}_M\cup\mathcal{V}_F$ and $\bm\alpha=\text{diag}(\alpha_i)$. 
\par
\begin{figure*}
\centering
\includegraphics[scale=1.13]{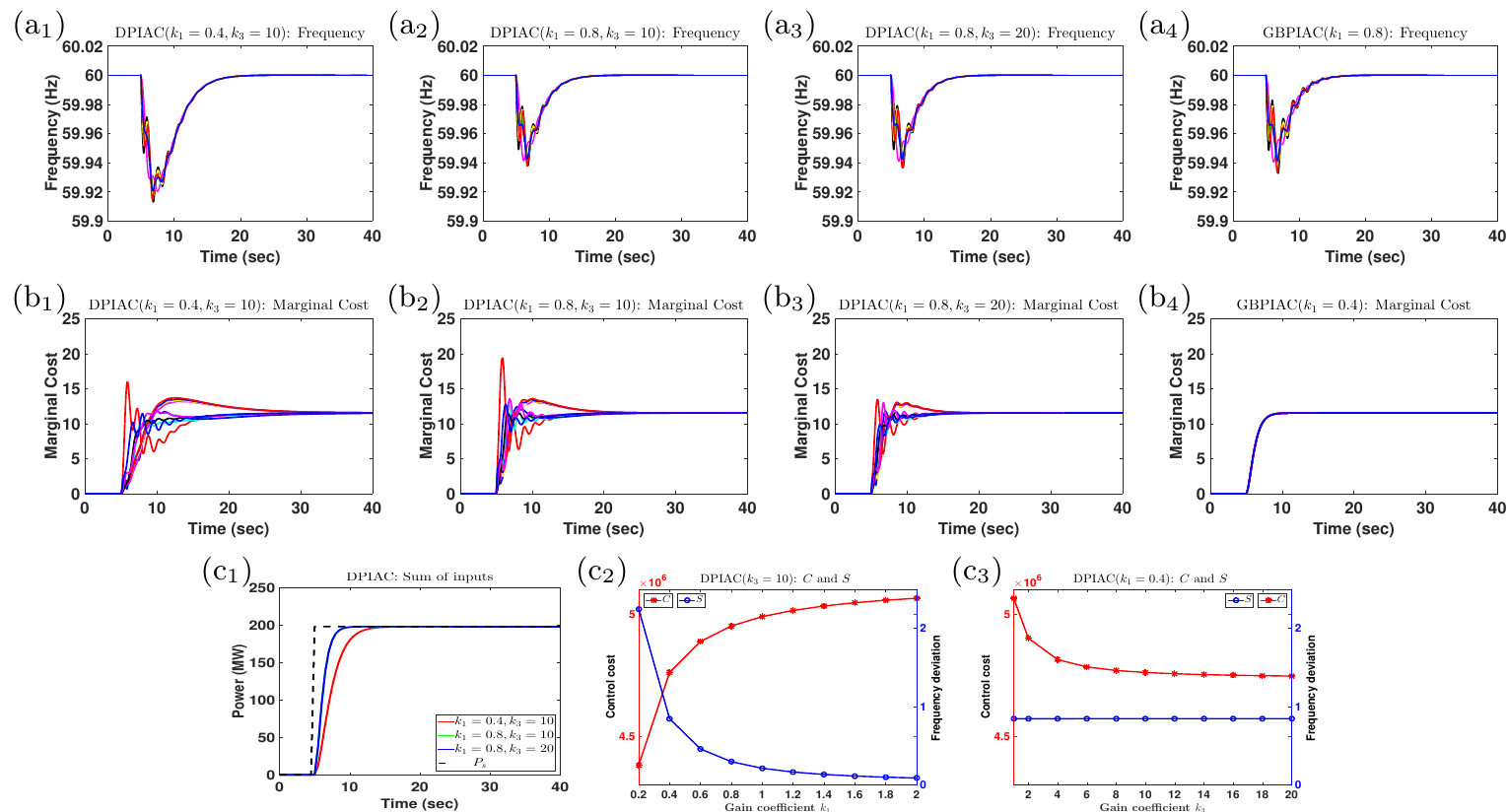}
\caption{The simulation result of the determistic system.}\label{fig1}
\end{figure*}

 
From Fig.\ref{fig1}~($a_1$-$a_2$), it can be observed that
the frequency restoration is accelerated by a larger $k_1$
with an accelerated convergence of the control input as shown in Fig. \ref{fig1}~($c_1$). From Fig.\ref{fig1}~($b_2$-$b_3$), it can be seen that the consensus of the marginal costs is accelerated by a larger $k_3$ with little influences on the frequency deviation as shown in Fig.\ref{fig1}~($a_2$-$a_3$).
It can be easily imagined that the marginal costs converge to that of GBPIAC as shown in Fig. \ref{fig1} (b4) as $k_3$ further increases. \emph{Hence, by increasing $k_1$ and $k_3$, the convergence of the state of the closed-sytstem to the optimal state can be accelerated, and by increasing $k_3$, the performance of the distributed control method DPIAC approaches to that of the centralized control method GBPIAC.} 
\par 
Fig.\ref{fig1}~($c_2-c_3$) show the trends of $S$ and $C$ with respect to $k_1$ and $k_3$. It can be observed from Fig.\ref{fig1}~($c_2$) that as $k_1$ increases,
the frequency deviation decreases while the control cost increases. \emph{Hence, to obtain a better performance of the frequencies,
a higher control cost is needed.} In addition, $S$ converges to a non-zero value as $k_1$ increases which is consistent with the anlsysis in (\ref{k1infinity}). However, the control cost is bounded as $k_1$ increases due to the bounded disturbance, which is different from the conclusion from Theorem \ref{theorem2}. When the disturbance is unbounded, the control cost is also unbounded, which will be further discussed in the next subsection. From Fig.\ref{fig1}~($c_3$), it can seen
that the control cost decreases inversely to a non-zero value as $k_3$ increases, which is also consistent with the analysis in (\ref{k3infinity}).
\begin{figure*}
\centering
\includegraphics[scale=1.13]{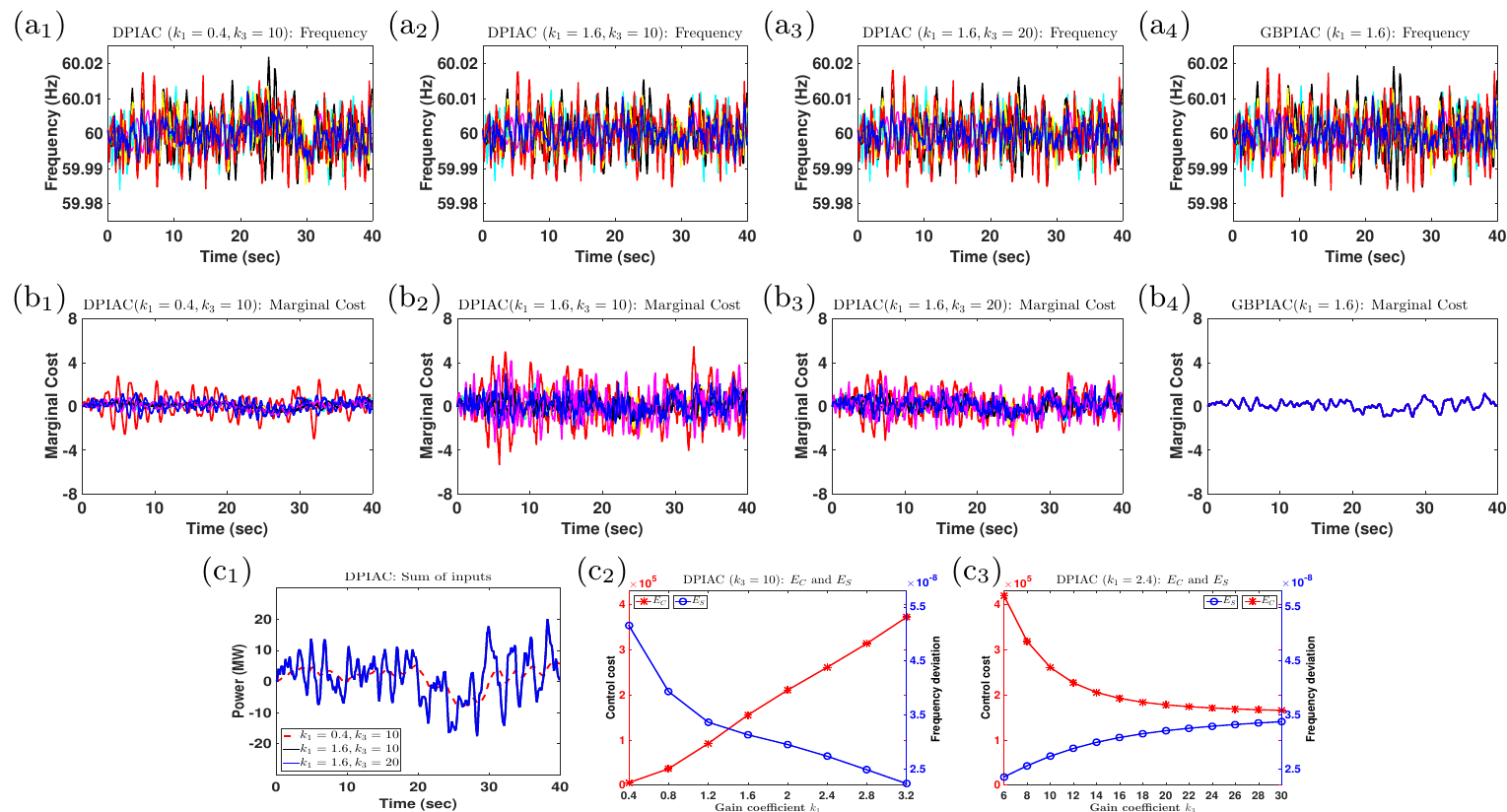}
\caption{The simulation result of the stochastic system.}\label{fig2}
\end{figure*}

\subsection{In the stochastic system}
With the interpretation of the $\mathcal{H}_2$ norm where the disturbances 
are modelled by white noise, we assume the disturbances are from the nodes of loads and $w_i\sim N(0,\sigma_i^2)$ with $\sigma_i=0.002$
for $i\in\mathcal{V}_P$. With these 
noise signal, the disturbances are unbounded and the closed-loop system becomes a stochastic algebraic differential system. We refer to \cite{algorithmStochastic} for a numerical algorithm to solve this stochastic system. 
\par 
It can be observed from Fig. \ref{fig2}~($a_1$-$a_2$) that the variance 
of the frequency deviation can be suppressed with a large $k_1$ which however
increases the variances of the marginal costs and the total control cost as shown in Fig. \ref{fig2}~($b_2$) and ~($c_1$) respectively. These observations are consistent with the analysis of Theorem \ref{theorem2} when $k_3$ is fixed.
From Fig. \ref{fig2}~($a_2$-$a_3$), it can be seen that increasing $k_3$ can effectively 
suppress the variances of the marginal costs. 
\par 
 We calculate the following 
metrics to study the impact of $k_1$ and $k_3$ on the variance of the frequency 
deviation and the expected control cost,
\begin{subequations}
\begin{align*}
E_S=E[\bm \omega^T(t)\bm\omega(t)],~\text{and}~E_C=\frac{1}{2}E[\bm u^T(t)\bm\alpha\bm u(t)].
\end{align*}
\end{subequations}
\par 
Fig.\ref{fig2}~($c_2$-$c_3$) show the trend of $E_S$ and $E_C$ as $k_1$ and $k_3$ increase. Similar to the discussion in the previous subsection, a trade-off can be found between the frequency deviation and the control cost in Fig.\ref{fig2}~($c_2$). The difference is that the control cost increases linearly as $k_1$ increases unboundly because of the unbounded disturbances. This is consistent with the result in Theorem \ref{theorem2}, which further confirms that a better
frequency response requires a higher control cost. 
\par 
Fig. \ref{fig2}~($c_3$) illustrates the trend of the expected control cost and  the variance of the frequency deviations with respect to $k_3$. It can be observed that the expected control cost decreases as $k_3$ increases, which is consistent as in (\ref{k3infinity}). However, the variance of the frequency deviation is slightly increased, which is also consistent with our analysis in (\ref{k1infinity}). 
\par 

\section{Conclusion}\label{chapter5:Sec:conclusion}

For the power system controlled by DPIAC, it has been demonstrated analytically and numerically that the transient performance of the frequency can be improved by tuning the coefficients monotonically, and a trade-off between
the control cost and frequency deviations has to be resolved to obtain a 
desired frequency response with acceptable control cost.

%
%

There usually are noises and delays in the state measurement and communications in practice, which are neglected 
in this paper. 
How these factors influence the transient behaviors of the state requires further investigation.



  \bibliographystyle{elsarticle-num-names} 
\bibliography{ifacconf} 
%

%
%
%
\end{document}